\documentclass[a4paper,reqno]{amsart}

\usepackage{amsmath,amsfonts}

\newtheorem{theorem}{Theorem}

\newtheorem{definition}{Definition}
\newtheorem{corollary}{Corollary}
\newtheorem{proposition}{Proposition}

\newcommand{\rest}{|}

\def \A{{\mathbb A}}
\def \B{{\mathbb B}}
\def \Z{{\mathbb Z}}
\def \E{{\mathbb E}}
\def \F{{\mathbb F}}
\def \R{{\mathbb R}}

\def \sphinf{S^{2}_{\infty}}
\def \Phiinf{\Phi_{\infty}}
\def \Einf{E_{\infty}}
\def \Ainf{A_{\infty}}
\def \threeball{\overline{B}^3}
\def \EY{\E\rest{\p X}}

\def \Ltwo{L^{2}}

\def \curlyA{{\mathcal{A}}}

\def \curlyC{{\mathcal{C}}}

\def \muo{\mu_{0}}
\def \ko{c_2(\E,f)[X]}

\def \perflat{2\pi / \muo}

\def \Iper{[0,\perflat]}

\def \rthree{{\mathbb R}^{3}}
\def \Xo{X^{o}}

\def \Daplus{D_{\A}^{+}}
\def \Dbplus{D_{\B}^{+}}
\def \Daminus{D_{\A}^{-}}
\def \Dk{D_{k}}
\def \Splus{S^{+}}
\def \Sminus{S^{-}}
\def \spinthree{S_{(3)}}

\def \trace{{\textrm {tr\ }}}
\def \endo{{\textrm {End}}}
\def \deg{{\textrm {deg~}}}

\def \ch{{\textrm{ch}}}
\def \dim{{\textrm{dim\ }}}
\def \ind{{\textrm{ind\ }}}
\def \supp{{\textrm{Supp}}}
\def \ker{{\textrm{ker\ }}}

\newcommand{\p}{\partial}

\begin{document}

\title[An $L^2$-Index Theorem for Dirac Operators on 
$S^{1}\times{\mathbb R}^3$]{An $L^2$-Index Theorem for Dirac Operators on 
$S^{1}\times{\mathbb R}^3$}

\author[T. M. W. Nye and M. A. Singer]{Tom M. W. Nye and Michael A. Singer}
\address{Department of Mathematics and Statistics\\ University of Edinburgh\\ James Clerk Maxwell Building\\ The King's Buildings\\ Edinburgh. EH9 3JZ.}

\email{tmwn@maths.ed.ac.uk \and michael@maths.ed.ac.uk}

%%%%%%%%%%%%

\begin{abstract}
An expression is found for the $L^2$-index of a Dirac operator 
coupled to a connection on a $U_n$ vector bundle over 
$S^{1}\times{\mathbb R}^3$. 
Boundary conditions for the connection are given which ensure the 
coupled Dirac operator is Fredholm.
Callias' index theorem is used to calculate the index when the 
connection is independent of the coordinate on $S^1$.
An excision theorem due to Gromov, Lawson, and Anghel reduces the 
index theorem to this special case. 
The index formula can be expressed using the adiabatic limit of the $\eta$-invariant of a Dirac operator canonically associated to the boundary. 
An application of the theorem is to count the zero modes of the Dirac operator in the background of a caloron (periodic instanton). 
\end{abstract}


\maketitle

\section{Introduction}

Let $X$ be a compact, oriented smooth manifold with boundary $\p
X$, and let $\Xo = X\backslash \p X$ be the corresponding open
manifold. 
Let $g$ be a complete Riemannian metric on $\Xo$ and let
$E\to \Xo$ be a complex vector bundle with a hermitian structure and
unitary connection $A$. If $X$ is a spin-manifold, we can introduce
the coupled Dirac operator 
$$D_A: C^\infty(\Xo,S\otimes E) \to C^\infty(\Xo,S\otimes E)
$$
and it is natural to attempt to obtain first,  conditions on $A$ and
$g$ near $\p X$ which ensure that $D_A$ is a Fredholm operator on
$\Ltwo$, 
and second, to obtain a formula for the $\Ltwo$-index.  Since with
standard conventions $D_A$ is self-adjoint, we must explain how
one obtains interesting $\Ltwo$-index problems from it.

If $\dim X$ is {\em even}, then the total spin-bundle decomposes as $S =
\Splus\oplus\Sminus$ and so $D_A$ gives rise to the `chiral' Dirac 
operators 
$$
D^{\pm}_A:C^\infty(\Xo,S^\pm\otimes E) \to C^\infty(\Xo,S^\mp\otimes 
E)
$$
which have equal and opposite $\Ltwo$-indices. If the geometry near 
$\p
X$ is restricted so that $g$ is a {\em $b$-metric}, then the
celebrated index theorem of Atiyah, Patodi and Singer \cite{ati73}  
gives a formula
for the $\Ltwo$-index of $D_A^+$ when this operator is Fredholm. (The
notion of $b$-metric was not used explicitly in \cite{ati73}; they 
worked with the equivalent idea of $X^o$ having cylindrical ends, a 
simple condition on $g$ near $\p X$.  The APS index theorem is 
discussed from the point of view of $b$-metrics in \cite{mel93}.) 
The APS theorem 
expresses the $\Ltwo$-index of $D_A^+$ as a sum of two terms, one an
integral of characteristic classes over $X$ and a boundary
contribution involving the famous $\eta$-invariant of $\p X$. 

By contrast, when $\dim X$ is {\em odd}, an interesting index-problem
arises for operators of the form
\begin{equation}\label{caltype}
D_{A,\Phi}:=D_A + 1\otimes\Phi: C^\infty(\Xo,S\otimes E) \to
C^\infty(\Xo,S\otimes E)
\end{equation}
where $\Phi$ is a suitable skew-adjoint endomorphism of $E$. 
According 
to work of Callias, Anghel and R{\aa}de \cite{cal78} \cite{ang93b} 
\cite{rad94}, $D_{A,\Phi}$ is Fredholm in
$\Ltwo$ under mild conditions on $\Phi$, the most important being that
it be invertible on $\p X$. Then, with no further restrictions on the
geometry of $g$ near $\p X$, the index is given by integrating over
$\p X$ a certain characteristic class constructed from $\Phi\rest{\p X}$ and
$E\otimes S\rest{\p X}$. We shall refer to this as the CAR index theorem 
in this paper, even though the result is closely related to 
pre-existing index theorems (cf.\ \cite{rad94}  for 
a discussion of this point).

The purpose of this note is to state and prove an index formula for
$D_A^+$ when $X$ is even-dimensional, but the geometry near the
boundary is not that of a  $b$-metric.  We shall restrict
ourselves to a very simple special case: we take $(\Xo,g)$ to be
isometric to $S^1\times \rthree$, with a standard flat product metric. 
Then 
$X$ is diffeomorphic to $S^1\times \threeball$ where $\threeball$ is the closed unit ball in $\rthree$.
Despite its
simplicity, this example already leads to an interesting index
theorem, thereby answering a question posed by Mazzeo and Melrose 
in 
their study of $\Phi$-pseudodifferential operators \cite{maz99}, at 
least in this very special case. 
(This $\Phi$ stands for `fibred-cusp' and has nothing to do with the $\Phi$ in 
\eqref{caltype}. This notational clash is unfortunate but seems 
unlikely to lead to serious confusion.)
We should note also that an $\Ltwo$-index theorem for the coupled Dirac
operator over $S^1\times S^1\times \R^2$ has been obtained by
Jardim \cite{jar99a}.  This is another natural example of a 
$\Phi$-geometry, but with fibre dimension $2$ rather than $1$.

The index theorem on $S^1\times\rthree$ is important in 
the study of self-dual calorons (periodic instantons) of which more 
will be said in \S\ref{sec:Nahm}.

\subsection{Statement of results}
\label{STR}
The following notation will be used throughout this paper:
let $X = S^1\times
\threeball$, $S^2_\infty = \p \threeball$, so that 
$\p X = S^1\times \sphinf$; let $p:X
\to \threeball$ be the projection. Let the metric $g$ on $\Xo$ be the
standard flat product metric on $S^1\times \rthree$ giving the circle 
a length $2\pi/\muo$, where $\muo >0$.
Let $z\in \R/ (2\pi/\muo)\Z$ be a standard  coordinate on the circle, 
and $x_1,x_2,x_3$ standard coordinates on $\rthree$.
Finally, fix an orientation on $X$ by decreeing that the ordered 
basis $(dz,dx_1,dx_2,dx_3)$ be positive.

Let $\E \to X$ be a smooth $U_n$-bundle and let $\A$ be a smooth
unitary connection in $\E$. $\A$ will be identified with the
corresponding covariant derivative operator $\nabla$, which has
components $\nabla_z,\nabla_1,\nabla_2,\nabla_3$ in the frame
$(\p_z,\p_1,\p_2,\p_3)$. Each of the two spin-bundles $S^\pm$ over
$X^o$ can be identified with $p^*S_{(3)}$, where $S_{(3)}$ is the
spin-bundle of $\rthree$. This is a complex vector bundle of rank 2
and comes equipped with skew-adjoint Clifford multiplication operators
$e_1,e_2,e_3$ associated with $\p_1,\p_2,\p_3$. The two coupled Dirac
operators over $X^o$ can now be written
\begin{equation}\label{eq:decompdirac}
D_{\A}^{\pm}= \pm\nabla_{z} + D_A
:C^\infty(\Xo, p^*S_{(3)}\otimes \E) \to C^\infty(\Xo,
p^*S_{(3)}\otimes \E)
\end{equation}
where $D_A = \sum_1^3 e_j\nabla_j$.  The first term in
\eqref{eq:decompdirac} operates on sections of the tensor product by
$\nabla_z(p^*(s)\otimes u) = p^*(s) \otimes \nabla_z u$ for any
$s \in C^\infty(\rthree,S_{(3)})$ and $u\in C^\infty(\Xo,\E)$.

Note that the  `3+1' decomposition of
$D_{\A}^{\pm}$ in \eqref{eq:decompdirac} 
depends only upon the product structure of $X$; we
introduced bases only for ease of presentation. 
 We now focus on
$\Daplus$; by abuse of notation, denote by the same symbol the
extension of $\Daplus$ to Sobolev spaces
\begin{equation} \label{sobop}
\Daplus: W^1(X,\Splus\otimes \E) \to W^0(X,\Sminus\otimes \E),
\end{equation}
where $W^k$ is the space of sections with $k$ derivatives in $\Ltwo$,
the latter space being defined in terms of the metric $g$.

In order to fix 
boundary conditions that make 
\eqref{sobop} a Fredholm operator it is convenient to fix a 
trivialization of $\E\rest{\p X}$:
\begin{definition}\label{def:framedbdl} A {\em framing} of  $\E$ is 
a choice of trivialization $f$ of $\EY$. The pair $(\E,f)$ is called a
{\em framed} bundle.
\end{definition}
As we shall see in \S\ref{k0deg}, framed bundles of rank $\geq 2$ are
classified by an integral topological invariant analogous to the
second Chern class, which we shall denote by $c_2(\E,f)[X]$.

For later convenience, we identify the trivial bundle that is 
implicit in Definition~\ref{def:framedbdl}
with $p^*\Einf$, where $\Einf$ denotes the trivial bundle
over $\sphinf$. Now we can write down boundary conditions for $\A$:

\begin{definition}\label{def:connection}
Let $\A$ be a connection on a framed bundle
$(\E,f)$, smooth up to the boundary of $X$,  let $A_\infty$ be a 
$U_n$-connection on $\Einf$ and let $\Phiinf$ be a skew-adjoint 
endomorphism of $\Einf$.
\begin{quote}
(i)  $\A$ is called a {\em
caloron configuration framed by $(A_\infty,\Phiinf)$} if
$$
\A = p^*A_\infty + p^*\Phiinf dz
$$
on $\p X$, where $f$ has been used to identify $\EY$ with $p^*E_\infty$.

(ii) The pair $(A_\infty,\Phiinf)$ is called {\em admissible} if 
$\nabla_\infty\Phiinf =0$, where $\nabla_\infty$ is the covariant 
derivative operator induced by $A_\infty$ on $\endo(\Einf)$.
\end{quote}
\end{definition}
Our first main result asserts that the Fredholm properties of \eqref{sobop}
are entirely determined by $\Phiinf$:
\begin{theorem} \label{thm:fredholm}
Let $\A$ be a caloron configuration framed by an admissible pair 
%overful hbox problems
\newline
$(A_\infty,\Phiinf)$. 
Then the operator in \eqref{sobop} is Fredholm if and only
if $1 - \exp(2\pi\Phiinf/\muo)$ is invertible.
\end{theorem}
If $(A_\infty,\Phiinf)$ is admissible, the eigenvalues
$i\mu_1,\ldots,i\mu_n$ of $\Phiinf$ are constant. 
An equivalent formulation  of this theorem is thus that
\eqref{sobop} is {\em not} Fredholm if and only if 
there exist $j\not=0$ and an integer $N$ such that $\mu_j = N\muo$.

We now come to a statement of our $L^2$-index theorem:
\begin{theorem} \label{thm:index}
Suppose that \eqref{sobop} is Fredholm. Then
\begin{equation} \label{infrm}
\ind(\Daplus) = -c_2(\E,f)[X] - \sum_k c_1(E^+_{(k)})[\sphinf]
\end{equation}
where for each integer $k$, $E^+_{(k)}$ is the 
sub-bundle of $\Einf$ on which $k\muo - i\Phiinf$ is positive-definite.
\end{theorem}
It is clear that $E^+_{(k)}$ is either $0$ or $\Einf$ for all but  a 
finite number of integers $k$.  Since $\Einf$ is trivial, it follows 
that the sum on the RHS of \eqref{infrm} is finite.

\begin{corollary}\label{cor:eta}
The index formula~$\eqref{infrm}$ can be re-expressed as
\begin{equation*}
\ind(\Daplus) =\int_X \ch(\E) -\frac{1}{2}\overline{\eta}_{\textrm{lim}}
\end{equation*}
where $\ch$ denotes the Chern character, and $\overline{\eta}_{\textrm{lim}}$ is the reduced $\eta$-invariant of the Dirac operator on $\partial X$ coupled to the connection $p^*A_\infty + p^*\Phiinf dz$ on $p^*\Einf$.
\end{corollary}

The Corollary is included as an indication of how the index theorem might be generalised and for comparison with the APS theorem---understanding $\eta$-invariants is not necessary for the proof of Theorem~\ref{thm:index}. 
We leave discussion of $\eta$-invariants to Appendix~\ref{sec:eta}, in which we sketch a proof of the Corollary.

\subsection{A sketch of the proof}

For the proof of Theorem \ref{thm:fredholm}, we apply general results 
of earlier
authors. One approach is to check the conditions written down by
Anghel, who has given necessary and sufficient conditions for the
Dirac operator over a complete Riemannian manifold to be Fredholm in
$\Ltwo$. The alternative is to use the characterisation of Fredholm
operators in the calculus of $\Phi$-pseudodifferential operators
\cite{maz99}. It is worth remarking that this latter approach gives
necessary and sufficient conditions for any `natural' operator on
$S^1\times \rthree$ to be Fredholm, not only operators of Dirac type. 
The two proofs appear in \S\ref{sec:fredholm}.

The proof of Theorem \ref{thm:index} involves two main steps. The 
first is a
calculation of the index in the case that there is a trivialization 
of $\E$ such that $\A$ is independent of 
$z$.  By Fourier analysis in the $S^1$-variable, the index can be 
identified in this case with a sum of indices of CAR-operators of the 
form \eqref{caltype}.
By topological arguments which we begin in \S\ref{topcal}, this 
calculation gives the index for any $\A$ over a framed bundle with 
$c_2(\E,f)[X]=0$.  

The second step invokes an excision theorem for
operators of Dirac type due to Anghel \cite{ang93} and Gromov--Lawson 
\cite{gro83}. 
In our case, this result gives
$\ind(D^+_{\A}) -
\ind(D^+_{\B})= -c_2(\E,f)[X]$
if $\B$ is any caloron configuration agreeing with $\A$ near $\p X$ 
but living on a new framed bundle $(\F,f)$, with $c_2(\F,f)[X]=0$. 
Since 
we calculated $\ind(D^+_{\B})$ in the first step, that completes the 
proof of  Theorem \ref{thm:index}. The details appear in \S\ref{indexp}.

\subsection{Application: counting the zero modes of the Dirac operator in the background of a caloron.}\label{sec:Nahm}

Connections $\A$ of the type we have considered here are of interest
in gauge theory, especially when they are required to satisfy the
self-duality equations \cite{gar88}; 
the term {\em caloron} was introduced in this context by Nahm \cite{nah83}.
The \emph{Nahm transform} of an self-dual caloron $\A$ is (roughly speaking) a bundle over $\R$ whose fibre at $t\in\R$ is the cokernel of the Dirac operator coupled to 
\begin{equation}\label{eq:defAt}
\A_t=\A-itdz.
\end{equation}
When $\A$ is self-dual a Weitzenb\"ock formula shows that $D^+_{\A_t}$ is injective, so the dimension of the cokernel is given by minus the index. 
At the time of writing, the Nahm description has only
been proved in detail for calorons on framed bundles with $c_2(\E,f)=1$ and such that the eigenbundles of $\Phiinf$ are trivial \cite{kra98b}.  Details of the full transform will be
the subject of a future publication. 

We want to express the index formula~$\eqref{infrm}$ in a form which is more familiar to mathematical physicists, and which is easier to interpret from the point of view of the Nahm transform. 
Fix a framed caloron configuration $\A$ on a framed bundle $(\E,f)$, and define $\A_t$ by~$\eqref{eq:defAt}$. 
$\Phiinf$ decomposes $\Einf$ into line bundles $\Einf= E_1\oplus\cdots\oplus E_n$ such that $\Phiinf$ acts by $i\mu_j$ on $E_j$, where $\mu_j\in\R$. 
Define $k_0=c_2(\E,f)$ and $k_j =c_1(E_j)[\sphinf]$ for $j=1,\ldots,n$. 
We can write the eigenvalues as $\mu_j=N_j \mu_0 +\epsilon_j$ for each $j$, where $0\leq \epsilon_j<\mu_0$ and $N_j\in\Z$. 
Moreover, we can assume that the eigenvalues are ordered so that $\epsilon_n\leq\epsilon_{n-1}\leq\ldots\leq\epsilon_1$ (this ensures our notation matches that of Garland and Murray \cite{gar88}). 
With this ordering, define $m_j=k_0+k_1+\cdots k_j$ for $j=0,\ldots,n-1$, and define intervals $I_{j,N}$ by
\begin{equation*}
I_{j,N}=
\begin{cases}
(\epsilon_1+N\muo,\epsilon_n+(N+1)\muo) &\textrm{when\ }j=0\textrm{\ and\ }N\in\Z \\
(\epsilon_{j+1}+N\muo,\epsilon_j+N\muo) &\textrm{for\ }j=1,\ldots,n-1,\textrm{\ and\ }N\in\Z.
\end{cases}
\end{equation*} 
Then Theorem~\ref{thm:fredholm} implies that $D^+_{\A_t}$ is Fredholm iff $t\in\bigcup I_{j,N}$, and Theorem~\ref{thm:index} implies that the fibre of the Nahm transform at $t\in I_{j,N}$ has rank $m_j$. 
Thus the rank depends periodically on $t$ and jumps at the points $t=N\mu_0 +\mu_j$ for $N\in\Z$ and $j=1,\ldots,n$. 
This agrees with the form of the Nahm transform anticipated in \cite[Section $8$]{gar88}.

\subsection{Acknowledgements} 
A slightly shorter version of this paper is to appear is the Journal of Functional Analysis---we thank the referee for his comments. 
Our thanks also go to Richard Melrose for making us
aware of the calculus of $\Phi$-differential operators and for several
enlightening conversations; and to Michael Murray for suggesting the problem of the Nahm transform for calorons.

\section{On the topology of calorons}
\label{topcal}

In this section we define the invariant $c_2(\E,f)[X]$ of a framed bundle
over $X$ from two points of view, the first homotopy-theoretic, the
second a version of Chern--Weil theory. Proposition~\ref{caltomon} 
and formula \eqref{caloroncharge} will be used in the proof of 
Theorem~\ref{thm:index}.

\subsection{Topological classification of framed bundles}
\label{k0deg}

We start with a useful way of thinking of framed
bundles and calorons in terms of the `rectangle' 
$R= [0,\perflat]\times\threeball$. There
is a natural map $R \to X$ which will be used to identify objects
defined over $X$ with corresponding objects over $R$. By abuse of
notation we denote the second projection of $R$ by $p$ and we shall
denote by the same symbol $(\E,f)$ the pull-back to $R$ of a framed
bundle over $X$; similarly for caloron configurations $\A$.  In
particular, when a framed bundle $(\E,f)$ is transferred to $R$,
we obtain a bundle over $R$, framed over $[0,\perflat]\times \sphinf$, and
with a `clutching map' 
$$
\phi: \E|\{0\}\times \threeball \simeq \E|\{\perflat\}\times \threeball.
$$
Since $\Einf$ is trivial, we can regard it as the restriction to
$S^2_\infty$ of a trivial bundle $E\to\threeball$, say, and we can
extend the framing $f$ of $\E$ to a bundle isomorphism $F:\E \to
p^*E$ over $R$. In this way $\phi$ becomes a unitary endomorphism $c$
of $E$ which shall refer to as a clutching function for $\E$. Because
of the periodicity, $c$ then lies in the group
\begin{equation*}
\curlyC = \{\mbox{unitary automorphisms }c\mbox{ of }E:
c\rest{\sphinf} = 1\}.
\end{equation*}
Now $\pi_0(\curlyC)=\Z$, for any element $c$ extends to a continuous
map from the one-point compactification $S^3$ of $\threeball$ into
$U_n$ and $\pi_3(U_n)=\Z$. We {\em define} $c_2(\E,f)[X]=-\deg(c)$.
$c_2(\E,f)[X]$ is the obstruction to extending the framing $f$ to the interior of $\E$---an extension exists iff $c_2(\E,f)[X]=0$.

We now reintroduce calorons, continuing to work over $R$, with $\E=p^*E$. 
Because of this identification, we can take the `3+1' decomposition
\begin{equation}\label{eq:splitting}
\nabla_{\A} = \nabla_{A_{(z)}} + dz(\p_z + \Phi_{(z)})
\end{equation}
along $\{z\}\times \threeball$, where $A_{(z)}$ is a unitary 
connection on $E$
and $\Phi_{(z)}$ is a skew-adjoint
endomorphism of $E$. Thus we have obtained from $\A$ a {\em path} 
$(A_{(z)},\Phi_{(z)})$ in
\begin{multline*}
\curlyA = \{(A,\Phi): A\mbox{ is a }U_n\mbox{ connection on }E,
\\ \Phi\mbox{ is a skew-adjoint endomorphism of }E, 
(A,\Phi)\rest{\sphinf} =(A_\infty,\Phiinf)\}.
\end{multline*}
For periodicity, the end-points of this path must be related by the 
clutching function $c$:
\begin{equation}\label{eq:clutchA}
A_{(\perflat)} = c^{\ast}(A_{(0)}) = cA_{(0)}c^{-1} - dcc^{-1}
\end{equation}
and 
\begin{equation}\label{eq:clutchphi}
\Phi_{(\perflat)} = c^{\ast}(\Phi_{(0)})=c\Phi_{(0)}c^{-1}.
\end{equation}
In other words, $\A$ can be identified with a loop in the quotient 
space $\curlyA/\curlyC$. Conversely, any such loop gives rise to a 
caloron configuration $\A$ framed by $(A_\infty,\Phiinf)$,
subject to further matching conditions needed to ensure the smoothness 
of $\A$ on $X$ when the two edges $\{z=0\}$ and $\{z = \perflat\}$ are
glued together.

The simplest example of this correspondence is of course the case that 
the path in $\curlyA$ is constant so that $c$ is identically 
$1$ and $c_2(\E,f)[X]=0$. Then we say that 
$\A$ is the {\em pull-back of a monopole}.  Here is a sort of converse:

\begin{proposition} Let $\A$ be a framed caloron in a framed
bundle with $c_2(\E,f)[X]=0$. Then there is a deformation $\B$ of $\A$ 
(through framed caloron configurations), such that $\B$ is the 
pull-back of a monopole.
\label{caltomon}\end{proposition}

\begin{proof} Since $c_2(\E,f)[X]=0$, we can find a unitary automorphism $C$ 
of $\E$ 
over $R$ which is equal to $1$ on $\{0\}\times\threeball$ and
$[0,\perflat]\times S^2_\infty$, and equal to $c$ on
$\{\perflat\}\times\threeball$. Pulling $\A$ back by 
$C$, we reduce to the
case $c=1$, so that the caloron configuration is now a {\em loop} in
$\curlyA$. But this space is contractible, so the result follows. \end{proof}

\subsection{Chern-Weil theory for framed bundles}\label{sec:chernweil}

Another way to think about the invariant $c_2(\E,f)[X]$ is in terms of
the integral
\begin{equation}\label{chargeintegral}
\int_{X} \ch (\E ) = -\frac{1}{8\pi^{2}}\int_{X} \trace 
F_{\A}\wedge F_{\A},
\end{equation}
where $\A$ is some framed caloron configuration. If $X$ had no boundary,
this integral would give minus the second Chern class, but
here there are additional contributions from $\p X$.
This integral has been calculated 
by different means in \cite{gro81} and  \cite{gar88}. As in the
previous section, we work over $R$, so that $\E$ is identified with
$p^*E$, together with a clutching function $c$.
Using the familiar trick of writing
\begin{equation*}
\trace F_{\A}\wedge F_{\A} = d~\trace \{d\A\wedge\A 
+\frac{2}{3}\A\wedge\A\wedge\A \}
\end{equation*}
the integral~$\eqref{chargeintegral}$ becomes an integral over the 
boundary of the rectangle $\Iper\times\threeball$:
\begin{equation}\label{eq:chernsimons}
-\frac{1}{8\pi^{2}}\int_{X} \trace F_{\A}\wedge F_{\A} =
-\frac{1}{8\pi^{2}}\int_{\p (\Iper\times\threeball)} \trace 
\{d\A\wedge\A +\frac{2}{3}\A\wedge\A\wedge\A\}.
\end{equation}
Now regard $\A$ as a path $(A_{(z)},\Phi_{(z)})$ satisfying
\eqref{eq:clutchA} and \eqref{eq:clutchphi}.
Evaluating~$\eqref{eq:chernsimons}$ on $(\p 
\Iper)\times\threeball$ and using the clutching formulas gives
\begin{multline*}
-\frac{1}{8\pi^{2}}\int_{(\p \Iper)\times\threeball} \trace 
\{d\A\wedge\A +\frac{2}{3} \A\wedge\A\wedge\A\} = \\
-\frac{1}{24\pi^2}\int_{\threeball}\trace (dcc^{-1})^3 + 
\frac{1}{8\pi^2}
\int_{\threeball}d~\trace \{A(0)c^{-1}dc\}.
\end{multline*}
The first term is $\deg c=-\ko$, and the second vanishes because 
$c=1$ on $\sphinf$.
On the other piece of the boundary we obtain
\begin{multline*}
 -\frac{1}{8\pi^{2}}\int_{\Iper\times\sphinf} \trace \{d\A\wedge\A 
+\frac{2}{3}\A\wedge\A\wedge\A \}= \\
 -\frac{1}{8\pi^{2}}\int_{\Iper\times\sphinf} \trace \{
2F_{A}\wedge\Phi dz - dA\wedge\Phi dz + A\wedge 
d\Phi\wedge dz +
\p_zA\wedge A\wedge dz\}.
\end{multline*}
The final term vanishes because the restriction of $A$ to $\Iper\times\sphinf$ is
pulled back from $\sphinf$ so that $\p_z A=0$ there (condition (i) of
Definition~\ref{def:connection}). On the other hand,  the sum  
of the  middle two terms is exact, so does not contribute to the integral.
The condition that $A_\infty$ is compatible with $\Phiinf$ implies that 
$A_\infty$ decomposes as a direct sum of connections, one on each eigenbundle of $\Einf$.
Suppose $E_\mu$ is the eigenbundle of $\Phiinf$ with eigenvalue $i\mu$. 
Then $A_\infty = \bigoplus a_\mu$ where $a_\mu$ is a connection on $E_\mu$, and $F_A \rest\sphinf = \bigoplus f_\mu$ where $f_\mu$ is the curvature of $a_\mu$. 
Since the first Chern class of $E_\mu$ is given by
\begin{equation}\label{eq:chernclass}
c_1(E_\mu)[\sphinf] = \frac{i}{2\pi}\int_{\sphinf} \trace f_\mu,
\end{equation}
we have 
\begin{equation*}
-\frac{1}{8\pi^{2}}\int_{\Iper\times\sphinf} \trace
2F_{A}\wedge\Phi dz = - \frac{1}{\muo}\sum_{\mu}\mu 
c_{1}(E_{\mu})[\sphinf].
\end{equation*}
Putting the terms together, we arrive at the expression
\begin{equation} \label{caloroncharge}
\int_{X} \ch(\E) = -\ko - \frac{1}{\muo}\sum_{\mu}\mu 
c_{1}(E_{\mu})[\sphinf].
\end{equation}

\section{Boundary conditions versus asymptotics}
\label{bcv}
In previous work, calorons have been  studied exclusively as connections 
over $S^1\times\R^3$, with decay conditions imposed near $\infty$;
 the compact manifold $X$ was not used. The purpose of this 
section is to show how the boundary conditions that are implicit in 
Definition~\ref{def:connection} translate into the `BPS' decay 
conditions for calorons that were written down in \cite{gar88}.

In order to compare the 
asymptotic region of $S^1\times \rthree$ with a neighbourhood of the 
boundary of
$X$, choose
polar coordinates $r,y_1,y_2$ in $\rthree$, and continue to use $z$ 
as a 
coordinate in 
$S^1$. Thus  $r$ is the distance from the origin in $\rthree$ and
$y_1$, $y_2$ are some local angular coordinates on $\sphinf$. We
suppose $y_1$ and $y_2$ are chosen so that $g$ takes the form
$$
g = dr^2 + r^2(h_{1}dy_1^2 + h_{2}dy_2^2) + dz^2,
$$
for some positive locally-defined functions $h_1,h_2$. Local
coordinates near the boundary of $X$ will be $x = r^{-1}$, $y_1,y_2$
and $z$, so that $x$ becomes a boundary defining function: $x\geq 0$
on $X$, with equality only at $\p X$. Writing $g$ in terms of $x$, 
$$
g = \frac{dx^2}{x^4} + h_1\frac{dy_1^2}{x^2} + h_2\frac{dy_2^2}{x^2}
+ dz^2.
$$
Now denote the components of $\A$, in some gauge (trivialisation) 
that is smooth up to the boundary, by 
$$
\nabla_x =\p_x + A_x,\,\nabla_{y_j}= \p_{y_j} + A_{y_j}, \nabla_z=\p_z + \Phi.
$$

Performing a $3+1$-decomposition of $\A= (A_{(z)},\Phi_{(z)})$ as before, 
we have, near the boundary,
\begin{align*}
\|A\|^2 &= |A_x|^2|dx|^2_g + |A_{y_1}|^2|dy_1|^2_g  + |A_{y_2}|^2|dy_2|^2_g \\
& =  x^2(x^2|A_x|^2 + h_1^{-1}|A_{y_1}|^2  + h_2^{-1}|A_{y_2}|^2).
\end{align*}
Since $x = r^{-1}$, we see that $\|A\| = O(r^{-1})$ as $r\to \infty$, 
uniformly in the angular variables $(y_1,y_2)$.  This statement is not 
gauge invariant: a better formulation is that there exist preferred
gauges near $\infty$ in $X^o$ (namely those that extend smoothly to
$\p X$), such that the connection $1$-form satisfies $\|A\| =
O(r^{-1})$ in such a gauge. In such gauges we also have $\|\Phi\| = O(1)$.

Similarly, the assumptions in Definition~\ref{def:connection} imply
that $\nabla_{y_j}\Phi$ and $\p_z A_{y_j}$ are $O(x)$ as $x\to 0$, while
$\nabla_x\Phi$ and $\p_z A_x$ are $O(1)$. It follows that
$\|\nabla_A\Phi - \p_z A\| = O(r^{-2})$ as $r\to \infty$, a fact that
will be used in the next section.

\section{Proof of Theorem~1.1}\label{sec:fredholm}
In this section we give two proofs of Theorem \ref{thm:fredholm} 
about the Fredholm
properties of $\Daplus$. The first proof rests on a result of Anghel 
in
\cite{ang93}, the second on the general theory of 
$\Phi$-pseudodifferential operators.  Of course both methods give the
same answer, and indeed the key point is the same in each case.

\subsection{First proof}

Theorem ($2.1$) of \cite{ang93} gives conditions for $D_{\A}=\Daplus\oplus\Daminus$ to be 
Fredholm:
$D_{\A}$ is Fredholm if and only if there is a compact set $K\subset \Xo$ and a
constant $C >0$ such that 
\begin{equation*}
\| D_{\A} \psi \|_{\Ltwo} \geq C \|\psi\|_{\Ltwo} ,\ \textrm{when}\ 
\psi\in W^{1} (S\otimes\E )\ \textrm{and}\ \supp (\psi )\subset 
\Xo\setminus K.
\end{equation*}
If $D_{\A}$ is Fredholm then $\Daplus$ must be Fredholm too.
Now 
\begin{equation*}
D_{\A}^{\ast}D_{\A} = (\Daminus\Daplus)\oplus 
(\Daminus\Daplus)^{\ast}
\end{equation*}
so to obtain estimates on $\| D_{\A}\|_{\Ltwo}$ we consider the 
operator $\Daminus\Daplus$.
Using the notation and conventions of \S\ref{STR} we have from 
\eqref{eq:decompdirac},
\begin{equation*}
\Daminus\Daplus = -\nabla_z^2+ [D_A,\nabla_z] + D_A^2.
\end{equation*}
The third term here is clearly positive, and
the boundary conditions allow us to  estimate the other two as 
follows.

{\textbf{The first term.}}
Extend the framing $f$ to a neighbourhood of $\p X$; 
this gives a gauge near $\infty$ in which the `3+1' decomposition~$\eqref{eq:splitting}$ can be performed.
As the boundary $\p X$ is approached the eigenvalues of $\Phi$ converge 
to the eigenvalues of $\Phiinf$.
Using spherical polar coordinates on $\rthree$, 
let $i\lambda_{j}(r,y_1,y_2,z)$ be the eigenvalues of $\Phi$, and 
$i\mu_{j}$ be the eigenvalues of $\Phiinf$ ($j=1,\ldots,n$) such that 
$\lambda_{j}\rightarrow\mu_{j}$ as $r\rightarrow\infty$.
Let $\lambda(r,y_1,y_2,z)$ be the smallest element in 
$\{ |\lambda_{j}+k\muo | : j=1,\ldots,n;~k\in\Z \}$ and $\mu$ be the 
smallest element of the set 
$\{ |\mu_{j}+k\muo | : j=1,\ldots,n,~k\in\Z \}$.
The invertibility condition on $\Phiinf$ in the statement of the 
theorem implies that $\mu > 0$,
so there exists a compact set $K_{1}\subset \Xo$ such that $\lambda > 
\mu / 2$ on $\Xo\setminus K_{1}$.

Suppose $\psi\in W^{2}(S^{+}\otimes\E)$ and $\supp~\psi\subset 
\Xo\setminus K_{1}$.
Using the isomorphism $\Splus\cong p^{\ast}\spinthree$,
$\psi$ can be written as a Fourier series
\begin{equation*}
\psi = \sum_{k}\exp(ik\muo z) \phi_k
\end{equation*}
where $\phi_{k}$ is a section of 
$\spinthree\otimes E$. Let
\begin{equation*}
\psi^{(k)} = \exp(ik\muo z) \phi_k.
\end{equation*}
Then
\begin{equation*}
  \nabla_{z}\psi^{(k)} = (ik\muo  +\Phi)\psi^{(k)}
\end{equation*}
so
\begin{equation*}
(-\nabla_{z}\nabla_{z}\psi^{(k)},\psi^{(k)})
\geq\frac{1}{4}{\mu}^{2}{\|\psi^{(k)}\|}^{2}, \qquad
\textrm{on\ }\Xo\setminus K_{1}
\end{equation*}
as a pointwise estimate.
(Since $\psi^{(k)}\in W^{2}(S^{+}\otimes \E)$, $\psi^{(k)}$ is 
actually continuous so both sides of the inequality exist.)
Since the inequality is independent of $k$ it holds for general $\psi$ 
and we obtain
\begin{equation} \label{fredhbd_1}
\supp (\psi)\subset \Xo\setminus K_{1} \Rightarrow
{(-\nabla_{z}\nabla_{z}\psi,\psi)}_{\Ltwo}
\geq\frac{1}{4}{\mu}^{2}{\|\psi\|}^{2}_{\Ltwo}.
\end{equation}

{\textbf{The second term.}}
We have
\begin{equation*}
[D_A,\nabla_z] = \sum_j e_j[\nabla_j,\nabla_z]
= \sum_j 
\iota(\p_j)(\nabla_A 
\Phi -\p_z A)
\end{equation*}
where $\iota(\xi)$ denotes interior product with $\xi$.
But $\|\nabla_{A}\Phi-\p_{z} A\| \rightarrow 0$ as $r\rightarrow 
\infty$, 
so there exists a compact set $K_{2}\subset \Xo$ such that
\begin{equation} \label{fredhbd_2}
\supp (\psi)\subset \Xo\setminus K_{2} \Rightarrow
| {([D_A,\nabla_z]\psi ,\psi )}_{L^2} | \leq
\frac{1}{8}\mu^{2}{\|\psi\|}^2_{L^2}.
\end{equation}

Now let $K$ be a compact set containing $K_{1}$ and $K_{2}$. 
Combining~$\eqref{fredhbd_1}$ and~$\eqref{fredhbd_2}$ we obtain 
\begin{equation*}
  \supp (\psi)\subset \Xo\setminus K\Rightarrow
  (\Daminus\Daplus\psi,\psi)_{\Ltwo}\geq 
\frac{1}{8}\mu^{2}{\|\psi\|}^{2}_{\Ltwo}.
\end{equation*}
A similar bound is obtained for 
$\Daplus\Daminus={(\Daminus\Daplus)}^{\ast}$, and so we obtain the 
following bound for $D_{\A}$:
\begin{equation*}
\psi\in W^{2}(S\otimes \E),~\supp (\psi)\subset \Xo\setminus 
K\Rightarrow
\|D_{\A}\psi\|_{\Ltwo}\geq\frac{1}{\sqrt{8}}\mu\|\psi\|_{\Ltwo}.
\end{equation*}
By density, the inequality in fact holds for $\psi\in W^{1}(S\otimes
\E)$. This completes the verification of Anghel's criterion and 
gives a proof of the `if' part of Theorem~\ref{thm:fredholm}.  The
`only if' part can also be proved in this framework but this is
omitted.  This converse statement also follows at once from the
discussion of the next section.

\subsection{Second proof}
Recall the boundary-adapted coordinate system
$x,y_1,y_2,z$ introduced in \S\ref{bcv}, and let the components of 
$\nabla_{\A}$ in these coordinates be
$$
\nabla_x =\p_x + A_x,\,\nabla_{y_j}= \p_{y_j} + A_{y_j}, \nabla_z=\p_z + A_z.
$$
Relative to a suitable choice of basis for the spin-bundles, we have
then
$$
\Daplus = \nabla_z + e_1x\nabla_{y_1} + e_2x\nabla_{y_2} +e_3 x^2\nabla_x.
$$
Strictly speaking, we are making a choice of normal coordinates here; 
otherwise there will be additional zero-order terms coming from 
connection coefficients.  This is an example of a $\Phi$-differential 
operator in the sense of \cite{maz99}; more generally the algebra of
$\Phi$-differential operators on $X$ consists of all differential 
operators which take the form
\begin{equation}\label{phidiff}
P(x,y,z;x^2\p_x,x\p_y,\p_z),
\end{equation}
near $\p X$, where $P$ is
smooth in the first three variables and polynomial in the last three
variables.  In \cite{maz99} it is shown that such an operator is 
Fredholm in $\Ltwo$ if and only if it is {\em fully elliptic} in the following 
sense. 
First, \eqref{phidiff} must be elliptic in the usual sense over 
$X^o$. Secondly, the associated {\em indicial family} must be 
invertible on every fibre $p^{-1}(y)\subset \p X$. Given such a fibre, the
indicial family on $p^{-1}(y)$ is defined by picking a real number 
$\xi$ and a real cotangent vector 
$\eta\in T_y^*\sphinf$, and defining
$$
\widehat{P}_{(y,\eta,\xi)} = P(0,y,z;i\xi,i\eta,\p_z)
$$
as a differential operator on $p^{-1}(y)$. To say that the indicial 
family is invertible is to say that
$\widehat{P}_{(y,\eta,\xi)}$ is invertible (in any reasonable 
space of sections over $p^{-1}(y)$), for each choice of
$(y,\eta,\xi)$ as above.

Following this recipe for $\Daplus$, we
obtain
$$
\widehat{P}_{(y,\eta,\xi)} =
\nabla_z + i(\eta_1e_1 + \eta_2e_2 + \xi e_3).
$$
This operator in $C^\infty(S^1, p^*S_{(3)}\otimes E_\infty)$ is a sum of 
two terms $B+A$, where $A=i(\eta_1e_1 + \eta_2e_2 + \xi e_3)$
is self-adjoint, $B= \nabla_z$ is skew-adjoint and 
$[A,B] =0$. It follows by considering $(A+B)^*(A+B)$
that $(A+B)u=0$ if and only 
if $Au=0$ and $Bu=0$. Now $B$ has a non-trivial null-space  only if
one of the $\mu_j$ is an integral multiple 
of $\mu_0$. Hence under the assumption of Theorem~1, $A+B$ is 
injective. Similarly the adjoint $(A+B)^*=A-B$ is injective, so that 
the hypothesis of Theorem \ref{thm:fredholm} implies that the 
indicial family is invertible, and so $\Daplus$ is Fredholm in
$\Ltwo$.  Conversely, if the condition fails, then $B$ is not
invertible, and nor is $B+A$ when $\eta_j=0=\xi$. So in this case
$\Daplus$ is not fully elliptic and hence cannot be Fredholm in
$L^2$.  The proof of Theorem~\ref{thm:fredholm} is now complete.

\subsection{Remarks about the $L^2$-condition}

According to \cite[Proposition 9]{maz99}, elements of the null-space
of a fully elliptic $\Phi$-differential operator decay very rapidly at
the boundary.  More precisely, if $\Daplus$ is fully elliptic,
if $\Daplus\psi =0$,  and if for some real  $m$,
$x^m \psi \in L^2(X)$,  then $\psi \in C^\infty(X)$ and $\psi$ 
vanishes to all orders in $x$ at $\p X$.
There is a similar statement for the cokernel. Now in terms of the 
boundary-adapted coordinates $(x,y_j,z)$, the volume element 
determined by the metric $g$ has the form $x^{-4}d\mu$
where $d\mu=h_1h_2dxdy_1dy_2dz$. It follows from the above that the 
index of \eqref{sobop} is the same as the index of
\begin{equation}\label{eq:ltwospaces}
\Daplus: W^1(X,\E\otimes S^+,d\mu) \to W^0(X,\E\otimes S^-,d\mu).
\end{equation}
This fact makes the next result almost obvious:
\begin{proposition}
Let $\A,\B$ be two caloron configurations on $(\E,f)$, both 
framed by $(A_\infty,\Phi_\infty)$. Then $\Daplus$ is Fredholm 
if and only if $D_{\B}^+$ is so, and their $L^2$-indices coincide.
\label{ddirac}    \end{proposition}
\begin{proof} The space of calorons on a given framed bundle, with given 
boundary data $(A_\infty,\Phi_\infty)$ is contractible. It is easy to 
see that any continuous path joining $\A$ to $\B$ gives rise to a 
norm-continuous path of Dirac operators between the Sobolev spaces in~$\eqref{eq:ltwospaces}$. 
Since each of these is 
Fredholm by Theorem~\ref{thm:fredholm}, it follows that the index is 
constant on this path. \end{proof}

\section{Proof of the index theorem}
\label{indexp}
\subsection{Proof when $\ko =0$}\label{proofkozero}

In this case, by  Propositions~\ref{caltomon} and \ref{ddirac}
it is enough to 
compute the index when $\E = p^*(E)$ and $\A = p^*A + p^*\Phi dz$ is 
the pull-back of a monopole (cf. \S\ref{k0deg}). Then the 
coefficients of $\Daplus$ are independent of $z$ and we can use
Fourier analysis in the 
$S^1$-variable to reduce the calculation of the index to that of a 
collection of operators of the form \eqref{caltype} on $\rthree$.
These operators are precisely the 
subject of the CAR theorem in its simplest form \cite{cal78}.

Let 
\begin{equation*}
  Y_{k}=\{\psi=\exp (ik\muo z)\phi:\phi\in W^0(\R^3,S_{(3)}\otimes E) \}
\end{equation*}
so that 
\begin{equation*}
  W^{0}(S^{+}\otimes \E) = \{ 
  \sum\psi^{(k)} : \psi^{(k)}\in Y_{k}\ \textrm{and\ } 
  \sum {\| \psi^{(k)} \|}^{2} < \infty \}.
\end{equation*}
Since by assumption the coefficients are independent of 
$z$, $\Daplus$ maps $Y_k\cap W^1$ into $Y_k$ and its restriction to 
this subspace is equal to
\begin{gather*}
\Dk : W^{1}(\spinthree\otimes E)\longrightarrow 
W^{0}(\spinthree\otimes E)
\\ \Dk =  D_{A}+ik\muo  + 1\otimes\Phi.
\end{gather*}
According to the theory developed by Callias--Anghel--R\aa de, 
$\Dk$ is Fredholm for every $k\in\Z$ iff $\Daplus$ is Fredholm, and
\cite{rad94} shows that:
\begin{align*}
\ind \Dk & = -\int_{\sphinf}\hat{A}(\sphinf)\wedge \ch(E_{(k)}^{+}) \\
& = -c_{1}(E^{+}_{(k)})[\sphinf]
\end{align*}
where $E^{+}_{(k)}$ is the subbundle of $\Einf$ on which 
$(k\muo-i\Phiinf)$ has positive eigenvalues.  (We have already noted 
that this sum is finite.)  Since $Y_j\cap Y_k = 0$ if 
$j\not=k$, the index of $\Daplus$ is the sum of the indices of the 
$\Dk$, i.e.
\begin{equation*}
\ind\Daplus=\sum_{k}\ind\Dk=-\sum_k c_{1}(E^{+}_{(k)})[\sphinf].
\end{equation*}
That completes the proof of Theorem~\ref{thm:index} when $c_2(\E,f)[X]=0$.

\subsection{Proof of the index theorem when $c_2(\E,f)[X]\not=0$}

Anghel \cite{ang93}, generalizing work of Gromov and Lawson \cite{gro83}, 
has given an excision theorem which compares the $L^2$-indices of a 
pair of Dirac operators over a complete manifold that agree near 
infinity. In our case this result yields the following statement. 
Let $\E$ and $\F$ be a pair of bundles over $X^o$ and let $\A$ and 
$\B$ be unitary connections on $\E$ and $\F$ (respectively). Suppose 
that there is a bundle isometry $\theta:\E|X^o\backslash K \to
\F|X^o\backslash K$ which carries $\A$ to $\B$ outside some compact set $K\subset\Xo$. 
Then 
\begin{equation}\label{Anghelsformula}
\ind\Daplus-\ind\Dbplus=\int_{X^o} \ch(\E)-\int_{X^o} \ch(\F).
\end{equation}

We are going to deduce Theorem~\ref{thm:index} by taking for $\B$ a 
connection which agrees with $\A$ near $\infty$, but which 
lives on a framed bundle $(\F,f)$ with $c_2(\F,f)=0$. This will 
complete the proof in view of the results of the previous section.

Let then $(\E,f)$ be a framed bundle and $\A$ a framed 
caloron configuration on $\E$. As in \S\ref{k0deg}, identify $\E$ with 
$p^*E$, together with a clutching function $c\in \curlyC$. Extend the 
framing smoothly from the boundary to a region $[0,\perflat]\times U$ where 
$K\subset \R^3$ is compact and $U=\Xo\setminus K$. 
By a deformation of $\A$ over 
$\{\perflat\}\times U$ which vanishes at $\infty$, we can assume that 
$c=1$ on $U$. Now define $\F=p^*E$ and $\B$ to agree with $\A$ 
over $S^1\times U$, but extended over $S^1\times K$ to define a 
smooth connection on $\F$. (This can be achieved by a suitable use 
of cut-off functions.)  

Applying \eqref{Anghelsformula},
\begin{equation*}
\ind\Daplus - \ind D_{\B}^{+} = \int_{\Xo}\ch (\E) - \int_{\Xo}\ch 
(\F).
\end{equation*}
But
\begin{equation*}
\int_{\Xo}\ch (\E)  =  -\ko - \frac{1}{\muo}\sum_{\mu}\mu
c_{1}(E_{\mu})[\sphinf]
\end{equation*}
from~$\eqref{caloroncharge}$, and
\begin{equation*}
\int_{\Xo}\ch (\F)  =  -\frac{1}{\muo}\sum_{\mu}\mu
c_{1}(E_{\mu})[\sphinf].
\end{equation*}
So
\begin{equation*}
\ind\Daplus = \ind D_{\B}^{+} - \ko
\end{equation*}
From \S\ref{proofkozero} we know that $\ind D_{\B}^{+} =-\sum_k
c_{1}(E^{+}_{(k)})[\sphinf]$ 
so we have proved that
\begin{equation*}
\ind\Daplus = -\ko -\sum_{k} c_{1}(E_{(k)}^{+})[\sphinf].
\end{equation*}
This completes the proof of Theorem \ref{thm:index}.

\appendix
\section{Adiabatic limits of $\eta$-invariants}\label{sec:eta}

The index formula~$\eqref{infrm}$ can be re-expressed in a form reminiscent of the APS formula by using adiabatic limits of $\eta$-invariants. 
The aim of this Appendix is to explain roughly how this can be done, thereby sketching a proof of Corollary~\ref{cor:eta}.

First recall the APS formula for the index of a Dirac operator $D$ on a manifold $X$ with cylindrical end \cite[Theorem $3.10$]{ati75}:
\begin{equation}\label{eq:APS}
\ind D = \int_X \alpha_0 - \frac{h+\eta(D_\p)}{2}.
\end{equation}
The first term is the integral over $X$ of some form $\alpha_0$, while the second depends on a Dirac operator $D_\p$ associated to the boundary $\p X$: 
$h=\dim\ker D_\p$ and $\eta(D_\p)$ is the $\eta$-invariant of $D_\p$. 
The $\eta$-invariant measures the asymmetry of the spectrum of $D_\p$, and is defined to be $\eta(D_\p)=\eta(D_\p)(0)$ where:
\begin{equation*}
\eta(D_\p)(s)=\sum_{\substack{\lambda\in\textrm{Spec~}D_\p \\ 
\lambda\neq 0}}(\textrm{sign}~\lambda)| \lambda |^{-s}.
\end{equation*}
Part of the proof of the APS theorem involves showing $\eta(D_\p)(s)$ is analytic at $s=0$, so that the definition $\eta(D_\p)=\eta(D_\p)(0)$ makes sense. 

Next consider a closed manifold $M$ equipped with a family of metrics $g_\epsilon$, where part of the metric blows up as $\epsilon\rightarrow 0$, and suppose that $D_\epsilon$ is a family of Dirac operators on $(M,g_\epsilon)$, with $\eta$-invariant $\eta(D_\p)$. 
Bismut and Cheeger \cite{bis89} show that, under certain conditions, not only does $\eta(D_\epsilon)$ exist for each $\epsilon>0$, but remarkably, the limit $\lim_{\epsilon\rightarrow 0}\eta(D_\epsilon)$ exists. 
They consider fibrations $Z\rightarrow M\xrightarrow{p} B$ of compact oriented spin manifolds and equip $M$ with a family of metrics
\begin{equation*}
g_\epsilon=\epsilon^{-1}p^*(g^B)+g^Z
\end{equation*}
where $g^B$ is a metric on $B$ and $g^Z$ annihilates the orthogonal complement of the fibres. 
When $D_\epsilon$ is a Dirac operator on $(M,g_\epsilon)$ coupled to some auxilliary bundle, Bismut and Cheeger \cite{bis89} show the limit of the reduced $\eta$-invariant, $\lim_{\epsilon\rightarrow 0}\overline{\eta}(D_\epsilon)$ exists. 
(The reduced $\eta$-invariant of an operator $D$ is $\overline{\eta}(D)=(h+\eta(D))/2$, which corresponds to the second term of the APS formula~$\eqref{eq:APS}$.)
Moreover, they show that the limit is given by the integral of some form $\hat{\eta}$ over the base $B$, and give explicit formulae for $\hat{\eta}$.

We apply this theory to the fibration $S^1\rightarrow S^1\times\sphinf\xrightarrow{p}\sphinf$, where the metric on $M=S^1\times\sphinf$ is given by 
\begin{equation*}
g_\epsilon=\epsilon^{-1}(h_1 dy_1^2 + h_2 dy_2^2) +dz^2
\end{equation*}
in the notation of Section~\ref{bcv}. 
Of course, this is the same as the restriction of the metric $g$ on $S^1\times\threeball$ to $\chi^2=\epsilon$. 
Let $D_\epsilon$ be the Dirac operator on $(M,g_\epsilon)$ coupled to the bundle $p^\ast\Einf$ via the connection $p^\ast\Ainf+p^\ast\Phiinf dz$. 
our aim is to calculate the limit $\overline{\eta}_{\textrm{lim}}=\lim_{\epsilon\rightarrow 0}\overline{\eta}(D_\epsilon)$. 
According to \cite[Theorem $4.95$]{bis89} this is given by an integral
\begin{equation*}
\overline{\eta}_{\textrm{lim}}=\frac{1}{2\pi i}\int_{\sphinf}\hat{\eta}. 
\end{equation*}
The definition of $\hat{\eta}$ \cite[Definition $4.93$]{bis89} simplifies in our situation to
\begin{equation}\label{eq:etahat}
\hat{\eta}=\frac{1}{\pi^{1/2}}\int_{0}^{\infty} 
{\textrm{tr}}^{\textrm{even}}~\big( D_Z\exp(-F_\infty-uD_Z^2) \big)\frac{du}{2u^{1/2}}
\end{equation}
where $F_\infty$ is the curvature of $\Ainf$ and $D_Z$ is the Dirac operator associated to the fibre $Z=S^1$, which is given by
\begin{equation*}
D_Z =-i\frac{d}{dz}-i\Phiinf.
\end{equation*}
Recall the notation of Section~\ref{sec:chernweil}: $\Phiinf$ decomposes $\Einf$ into a direct sum of eigen-bundles $\Einf=\bigoplus E_\mu$. 
The Fredholm condition (Theorem~\ref{thm:fredholm}) implies that each eigenvalue $i\mu$ of $\Phiinf$ can be written as $\mu=N_{\mu}\mu_0+\epsilon_\mu$ where $N_\mu\in\Z$ and $0<\epsilon_\mu<\mu_0$. 
$D_Z$ decomposes as $D_Z=\bigoplus D_Z^{(\mu)}$, where
\begin{equation*}
D_Z^{(\mu)} =-i\frac{d}{dz}+\mu
\end{equation*}
and $\Ainf$ decomposes as a direct sum of connections on each eigen-bundle, $\Ainf=\bigoplus a_\mu$. 
Then
\begin{align*}
\overline{\eta}_{\textrm{lim}}&=\frac{1}{2\pi i}
\int_{\sphinf}\frac{1}{\pi^{1/2}}\int_{0}^{\infty}
{\textrm{tr}}^{\textrm{even}}~\big( D_Z\exp(-F_\infty-uD_Z^2) \big)\frac{du}{2u^{1/2}} \\
&=\sum_{\mu}\bigg(\frac{i}{2\pi}\int_{\sphinf}f_\mu\bigg)
\bigg(\frac{1}{\pi^{1/2}}\int_0^\infty
\textrm{tr~}\big[D_Z^{(\mu)}\exp(-u(D_Z^{(\mu)})^2)\big]\frac{du}{2u^{1/2}} \bigg)
\end{align*}
where $f_\mu$ is the curvature of $a_\mu$. 
Using equation~$\eqref{eq:chernclass}$ the first bracket is $c_1(E_\mu)[\sphinf]$ while the second is the $\eta$-invariant of $D_Z^{(\mu)}$ which we denote $\eta_\mu$ (this follows from \cite[Equation $0.3$]{bis89} or \cite[Theorem $2.6$]{bis86}). 
Hence
\begin{equation*}
\overline{\eta}_{\textrm{lim}}=\sum_{\mu}\eta_\mu c_1(E_\mu)[\sphinf].
\end{equation*}
Calculating $\eta_\mu$ is a standard example (see \cite{ati73} or \cite[Section $1.10$]{gil95}):
\begin{equation*}
\eta_\mu = 1-\frac{2\epsilon_\mu}{\mu_0}
\end{equation*}
and so 
\begin{equation}\label{eq:adiabatic}
\overline{\eta}_{\textrm{lim}}=-\frac{2}{\mu_0}\sum_{\mu}\epsilon_\mu c_1(E_\mu)[\sphinf].
\end{equation}
since 
\begin{equation}\label{eq:trivsum}
\sum_\mu c_1(E_\mu)[\sphinf]=0.
\end{equation}

Using~$\eqref{caloroncharge}$, the index formula~$\eqref{infrm}$ can be written as
\begin{equation*}
\ind(\Daplus) = \int_X\ch(\E)+\frac{1}{\mu_0}\sum_\mu \mu c_1(E_\mu)[\sphinf] - \sum_k c_1(E^+_{(k)})[\sphinf].
\end{equation*}
The last two terms together become
\begin{multline*}
\frac{1}{\mu_0}\sum_\mu \mu c_1(E_\mu)[\sphinf] - \sum_k c_1(E^+_{(k)})[\sphinf]
 = \frac{1}{\mu_0}\sum_\mu \epsilon_\mu c_1(E_\mu)[\sphinf]\\
+\Big( \sum_{\mu}N_{\mu} c_1(E_\mu)[\sphinf]-
\sum_k \sum_{\mu : N_\mu\geq -k}c_1(E_\mu)[\sphinf] \Big).
\end{multline*}
The term in brackets vanishes because of identity~$\eqref{eq:trivsum}$, so
\begin{align}
\ind(\Daplus)&= \int_X\ch(\E) + \frac{1}{\mu_0}\sum_\mu \epsilon_\mu c_1(E_\mu)[\sphinf]\notag \\
&= \int_X\ch(\E)-\frac{1}{2}\overline{\eta}_{\textrm{lim}}\label{eq:final}
\end{align}
using equation~$\eqref{eq:adiabatic}$. 
Thus we have proved Corollary~\ref{cor:eta}. 
In this form the index formula resembles the APS formula~$\eqref{eq:APS}$ and it seems likely that equation~$\eqref{eq:final}$ might apply more widely to Dirac operators on manifolds with fibred boundaries. 

%% Bibliography
\bibliographystyle{amsplain}

\providecommand{\bysame}{\leavevmode\hbox to3em{\hrulefill}\thinspace}

\end{document}